\documentclass{amsart}

\usepackage{latexsym}
\usepackage{amsfonts}
\usepackage{amssymb}
\usepackage{color}
\usepackage{eufrak}
\newtheorem{teor}{Theorem}[section]
\newtheorem{lema}[teor]{Lemma}
\newtheorem{prop}[teor]{Proposition}

\newtheorem{defi}[teor]{Definition}
\newtheorem{afi}{Claim}

\newtheorem{ex}[teor]{Example}

\newcommand{\N}{\mathbb{N}}

\begin{document}

\title{Spaceability in sets of operators on $C(K)$}

\author[R. Fajardo]{Rog\'erio Fajardo}
\address{Instituto de matem\'atica e estat\'istica, Universidade de S\~ao Paulo, Rua do Mat\~ao, 1010, CEP 05508-900, S\~ao Paulo,  Brazil}
\email{leonardo@ime.usp.br}

\author[P. L. Kaufmann]{ Pedro L. Kaufmann }
\address{Instituto de matem\'atica e estat\'istica, Universidade de S\~ao Paulo, Rua do Mat\~ao, 1010, CEP 05508-900, S\~ao Paulo,  Brazil}
\email{plkaufmann@gmail.com}
\thanks{The second author was supported by CAPES, Research Grant PNPD 2256-2009.}

\author[L. Pellegrini]{ Leonardo Pellegrini}
\address{Instituto de matem\'atica e estat\'istica, Universidade de S\~ao Paulo, Rua do Mat\~ao, 1010, CEP 05508-900, S\~ao Paulo,  Brazil}
\email{leonardo@ime.usp.br}

\maketitle

\begin{abstract}  We show that when  $C(K)$ does not have few operator -- in the sense of Koszmider (\cite{Ko}) -- the sets of operators which are not weak multipliers is  spaceable. This shows a contrast with what happens in general Banach spaces that do not have few operators. 

In addition,  we show that there exist a  $C(K)$ space  such that each operator on it  is of the form $gI+hJ+S$, where $g,h\in C(K)$ and $S$ is strictly singular, in connection to a result by Ferenczi (\cite{Fe}).

\end{abstract}



\section{Introduction and terminology}

Banach spaces with few operators have been the subject of several papers. 
In \cite{GM}, Gowers and Maurey build an hereditarily indecomposable Banach space such that every (bounded) linear  operator on it has the form  $\lambda I +S$, where $\lambda\in  \mathbb{C}$ and $S$ is strictly singular. 
Recently, in \cite{AH}, Argyros and Haydon have built a Banach space in which every operator has the form $\lambda I +C$, where $\lambda\in  \mathbb{R}$ and $C$ is a compact operator.  

We  restrict ourselves  to Banach spaces of the form  $C(K)$,   the space of real continuous functions on a Hausdorff compact topological space $K$, normed by the supremum.

For any g in $C(K)$, the multiplication operator $gI : f \mapsto gf$ is a  bounded linear operator on  $C(K)$.   In addition, it can be shown (see \cite{Ko}) that  for all $C(K)$ there is always an operator  which cannot be expressed in the form $ gI + C$, where $g\in  C(K)$ and $C$ is compact and  there is always an operator  which cannot be expressed in the form $ \lambda I + S$, where $\lambda\in  \mathbb{R}$ and $S$ is strictly singular. 

On the other hand,  Koszmider proved in  \cite{Ko} that there exists a $C(K)$ space  on which every operator has the form $gI+S$, where $g\in C(K)$ and
$S$ is strictly singular or, equivalently, weakly compact (the equivalence holds for $C(K)$ spaces, see \cite{Pe}).  Such operators are often called {\it weak multiplications} and Banach spaces $C(K)$ which every operator on it is a   weak multiplications are said to have {\it few operators}. If K is connected and all operator on   $C(K)$ is a   weak multiplications then   $C(K)$ is  indecomposable (see \cite{Fa}). Note that, since every  $C(K)$  contains a copy of $c_0$, it is not hereditarily indecomposable.

In what concerns  $C(K)$ spaces with few operators, it is usual to consider  a slightly larger class of operators than weak multiplications: it  is usual to consider the  class of  the {\it weak  multipliers} (see Definition \ref{defi:multiplicadorfraco}). A Banach space $C(K)$ on which every operator is a   weak  multiplier  is also  commonly said to have {\it few operators}.
Every weak multiplication is a weak multiplier (see Proposition \ref{prop:multiplicationmultiplier}). However, the notion of weak multiplier seems to be more appropriate. For example, if every operator on $C(K)$ is a weak multiplier and $C(K)$ isomorphic to $C(L)$, then every operator on  $C(L)$ is also a weak multiplier;   by the other hand, any $ C(K)$ (having few operators or not) is isomorphic to some $ C(L)$ which admits a weak multiplier which is not a weak multiplication (see \cite{Sc}). Furthermore,  weak multipliers  play a central role in all constructions of $C(K)$'s with few operators, as in \cite{Ko} and \cite{Pl}. Throughout this work, by \emph{$C(K)$ has few operators} we will mean that \emph{each operator in $C(K)$ is a multiplier}.


The question that motivated this work is the following: when $C(K)$ does not have few operators, \emph{how many operators which are not weak  multiplier exist in $C(K)$?}

In order to measure the amount of non weak  multipliers on a $C(K)$ space 
we use the concepts of lineability and spaceability, which appears in several papers: a subset $M$ of a Banach space $X$ is said to be \emph{lineable} if $M\cup \{0\}$ contains an infinite dimensional subspace of $X$. If that subspace can be chosen with  the additional property of being closed, then $M$ is said to be \emph{spaceable}. This terminology was firstly published in \cite{ags}, but it first appeared  in unpublished notes by Enflo and Gurariy, which will now formally appear as a part of \cite{egs}, a work in collaboration with Seoane-Sep\'ulveda.

Among other results, we  will show that when $C(K)$ does not have few operators, then the set of the operators on $C(K)$ which are not  weak  multiplier is spaceable -- or equivalently, that the quotient of the space of all operators on $C(K)$ by the subspace of weak    multipliers is infinite-dimensional. This equivalence is shown in Theorem \ref{teor:equivalent}.  Our results reveal a contrast between the few operators theory in general Banach spaces and in $C(K)$.   In a general Banach space setting, Ferenczi \cite{Fe} provides an example where the quotient of  the space of all operators on an Banach space $X$ by the subspace of the operators of the form
$\lambda I+S$, where $\lambda\in\mathbb{R}$ and $S$ is strictly singular,  has dimension 1. 
In other words, that there is an operator $J$ on $X$ such
that every operator on $X$ has the form $\lambda_1 J+\lambda_2I+S$,
where $\lambda_1,\lambda_2\in\mathbb{R}$ and $S$ is  strictly singular.
 We show that this does not occur for  $C(K)$ spaces.  Although, we get a related result if we allow replacing scalar numbers by continuous functions:
for any  $C(K)$ such that  every operator is a   weak  multiplication, there is an operator $J$ on $C(K\times\{0,1\})$ which is not weak multiplier such that every operator on $C(L)$ has the form  $gJ+hI+S$, where $g,h\in C(K\times\{0,1\})$ and $S$ is strictly singular.

All topological spaces referred in this paper are Hausdorff. We use  $\N$ to denote the set of non-negative integers and ${\mathcal L}(X)$ to denote the space of bounded linear operators on a Banach space $X$.

\section{Weak multipliers and weak multiplications}\label{sec:multipliers}

In this section we present the definitions of weak multipliers and weak multiplications, as well as some basic results on these kinds of operators that will be  used in this paper. To know more about how weak multipliers and weak multiplications are  related we mention \cite{Ko} and \cite{Sc}.

\begin{defi}\label{defi:multiplicacaofraca}\emph{\textbf{(\cite{Ko}, 2.1)}}\emph{
An operator $T:C(K)\longrightarrow C(K)$ is a weak multiplication if there is $g\in C(K)$
and an weakly compact operator $S$ on $C(K)$ such that $T=gI+S$.
}\end{defi}

We say that a sequence $(e_n)_{n\in\N}$ in $C(K)$ is 
\emph{pairwise disjoint} if $e_i\cdot e_j=0$, for all $i\neq j$.

\begin{defi}\label{defi:multiplicadorfraco}\emph{\textbf{(\cite{Ko}, 2.1)}}\emph{
An operator $T:C(K)\longrightarrow C(K)$ is a weak multiplier if,  
for every pairwise disjoint bounded sequence  $(e_n)_{n\in\N}$  in $C(K)$ 
and every sequence $(x_n)_{n\in\N}\subseteq K$ satisfying $e_n(x_n)=0$, we have
$$\lim_{n\rightarrow\infty}T(e_n)(x_n)=0.$$ 
}\end{defi}

Every weak multiplication is a weak multiplier. To show this, we use
the following caracterization of weakly compact operators on $C(K)$.

\begin{teor}\label{teor:weaklycompact}\emph{\textbf{(\cite{DU}, VI, Cor.17)}} 
Let $X$ be a Banach space. An operator $T:C(K)\longrightarrow X$  is weakly compact if and only if for every
bounded pairwise disjoint sequence $(e_n)_{n\in\N}$ we haver $||T(e_n)||$ converges to 0.
\end{teor}

\begin{prop}\label{prop:multiplicationmultiplier}
 Every weak multiplication is a weak multiplier.
\end{prop}
\paragraph{Proof:} 
Let $T=gI+S$ be a weak multiplication on $C(K)$, where $g\in C(K)$ and $S$ is a weakly compact operator. Let $(e_n)_{n\in\mathbb{N}}$ be a bounded pairwise disjoint sequence in $C(K)$ and let $(x_n)_{n\in\N}$ be a sequence in $K$. 
Then $T(e_n)(x_n)=g(x_n)\cdot e_n(x_n)+S(e_n)(x_n)=S(e_n)(x_n)$. 
By Theorem~\ref{teor:weaklycompact} we have that $\lim_{n\in\N}||S(e_n)||=0$.
Therefore $\lim_{n\rightarrow\infty}T(e_n)(x_n)=0$, proving that $T$ is a weak multiplier. 
\hfill{$\blacksquare$} 

\vspace{2mm}

We will also make use of the following:

\begin{teor}\label{teor:isosobre}\emph{\textbf{(\cite{Ko}, 2.3)}} 
If $K$ does not have a non-trivial convergent sequence and $T:C(K)\longrightarrow C(K)$
is a weak multiplier, then $T$ is onto $C(K)$ if and only if it is an isomorphism onto its range.
\end{teor}

Let ${\mathcal L}(C(K))$ be the Banach space of all bounded operators on $C(K)$, and let $M(C(K))$ be the subspace of ${\mathcal L}(C(K))$ consisting of the weak multipliers on $C(K)$. 

\begin{teor}\label{teor:closed}
$M(C(K))$ is closed in ${\mathcal L}(C(K))$. 
\end{teor}

\paragraph{Proof:} Let $(e_n)\in C(K)$ and $(x_n)\in K$ be sequences such that $M:=\hbox{sup}\|e_n\|<\infty$,  
$e_n\cdot e_m=0$ and $e_n(x_n)=0$. Take $S\in \overline{{\mathcal M}}$. 
Given $\varepsilon >0$, there exists $T\in {\mathcal M}$ satisfying $\| S-T\|<\frac{\varepsilon}{2M}$. 
Since $T\in  M(C(K))$ , there exists $N\in\mathbb{N}$ such that $|T(e_n)(x_n)|< \frac{\varepsilon}{2}$, for all $n>N$. Hence, for $n>N$, we have 
$$|S(e_n)(x_n)|\leq |T(e_n)(x_n)|+ |T(e_n)(x_n)- S(e_n)(x_n) | 
<\frac{\varepsilon}{2} + \| S-T\|\|e_n\| <\varepsilon, $$
which proves that $S(e_n)(x_n)\rightarrow 0$; therefore, $S\in  M(C(K))$. 
\hfill{$\blacksquare$}

\section{Lineability and spaceability of sets of non-weak multipliers}

We recall that a subset $S$ of a Banach space $X$ is \emph{lineable} if $S\cup\{0\}$ contains
an infinite-dimensional subspace of $X$. We say that $S$ is \emph{spaceable} if $S\cup\{0\}$ contains a closed infinite-dimensional subspace of $X$. In case $S$ is a complement of a closed subspace of $X$, the next theorem -- which is a consequence of a Theorem of \cite{KT} -- prove that both definitions are equivalent.

\begin{teor}\label{teor:equivalent} Let $M$ be a  subspace of a vector space $X$. 
The following assertions are equivalent:
\begin{itemize}
\item[(a)] $X\smallsetminus M$ is lineable;
\item[(b)] The quotient $X/M$ is infinite-dimensional;
\item[(c)] There are no $x_1,\ldots, x_n\in X$ such that for every 
		$x\in X$ there are
		$c_1,\ldots,c_n\in\mathbb{R}$ and $y\in M$ such that 
		$x=c_1x_1+\ldots+c_nx_n+y$.
\end{itemize}
If,  in addition,  $X$ is a Banach space and $M$ is closed, the assertions above are equivalent to

(d) $X\smallsetminus M$ is spaceable.
\end{teor} 

\paragraph{Proof:} The equivalence
between $(b)$ and $(c)$ follows immediately from the definition of quotient. We will prove that $(a)$ and $(b)$ are equivalent:

Let $\{z_n+M :n\in\mathbb{N}\}$ be a linearly independent family in $X/M$. We will prove that 
$\{z_n :n\in\mathbb{N}\}$ is linearly independent in $X$. In fact, 
\begin{eqnarray*}
\sum_{n=0}^k \alpha_nz_n=0 &\Rightarrow& \sum_{n=0}^k \alpha_nz_n\in M 
\Rightarrow\big(\sum_{n=0}^k \alpha_nz_n\big)+M = 0 \\
&\Rightarrow&  \sum_{n=0}^k \alpha_n\big(z_n+M\big) = 0 \Rightarrow \alpha_n=0, 
\forall n=0,1,\ldots , k.
\end{eqnarray*}
Therefore,  $S=\hbox{span}\{ z_n :n\in\mathbb{N}\}$ is infinite-dimensional. 
Let $\sum_{n=0}^k \alpha_nz_n$ be an element of $S$. 
If  $\sum_{n=0}^k \alpha_nz_n\in M$, then $ \sum_{n=0}^k \alpha_n\big(z_n+M\big) = 
\big( \sum_{n=0}^k \alpha_n z_n\big)+M =0$, 
and hence $\alpha_n=0$, for all $n=0,1,\ldots , k$. 
Therefore  $S\subset (X\smallsetminus M) \cup \{0\}$.

On the other hand, let  $Z$ be a  infinite dimensional subspace contained in  $(X\smallsetminus M) \cup \{0\}$,  take $\{z_n:n\in N\}$ a linearly independent subset of $Z$, and consider the set of classes $\{z_n+M :n\in\mathbb{N}\}\subset X/M$. 
Then if $\alpha_0, \alpha_2, \ldots, \alpha_k$ are scalars, we have that 
\begin{eqnarray*}
\sum_{n=0}^k\alpha_n(z_n+M)=0 &\Rightarrow& \sum_{n=0}^k (\alpha_nz_n)+M=0  \Rightarrow \sum_{n=0}^k \alpha_nz_n\in M \\
&\Rightarrow&  \sum_{n=0}^k \alpha_nz_n =0 \Rightarrow \alpha_n=0, \forall n=0,1,\ldots , k,
\end{eqnarray*}
and therefore $\{z_n+M :n\in\mathbb{N}\}$ is  linearly independent.

If  $X$ is a Banach space and $M$ is closed subspace, the equivalence between $(b)$ and $(d)$ follows from \cite[Theorem 2.2]{KT}.
\hfill{$\blacksquare$}

\vspace{2mm}

Let $NM(C(K))$ be the set of all operators on $C(K)$ which are not weak multipliers. If $C(K)$ has few operators then, by definition $NM(C(K))=\emptyset$. Otherwise, we will show that  $NM(C(K))$ is spaceable.  We start by studying a particular case.

\begin{teor}\label{teor:convergente} If $K$ has a non-trivial convergent sequence, then  $NM(C(K))$ is spaceable. 
\end{teor}

\paragraph{Proof:}

Let $(x_k)_{k\in \N}$ be a sequence of distinct elements of $K$ which converges to 
$\overline{x}$.

Using the normality of $K$, we get a sequence of continuous functions 
 $f_k:K\longrightarrow [0,1]$
with pairwise disjoint supports, such that $f_k(x_k)=1$ and $f_k(\overline{x}_l)=0$ 
for all $k\neq l$. We may assume that $f_k(x)=0$, for all $k$.

For each $n \geq 1$,  define the operator $T_n:C(K)\longrightarrow C(K)$ by
$$T_n(f)(x)=\sum_{j=0}^\infty f_j(x)(f(x_{j+n})-f(\overline{x})).$$
It is easy to verify that $T_n$ is well defined and $||T_n||=2$.
Let $T =\sum_{i=1}^{n} \lambda_iT_i$. For every $f$ in the unit ball of $C(K)$ 
and every $x\in K$, we have that
\begin{eqnarray*}|T(f)(x)|&=&|\sum_{j=0}^\infty \sum_{i=1}^n (\lambda_i f_j(x)(f(x_{j+i})-f(\overline{x})))|\\
&=&|\sum_{j=0}^\infty (f_j(x)\sum_{i=1}^n \lambda_i(f(x_{j+i})-f(\overline{x})))|\\
&=&|f_0(x) \sum_{i=1}^n \lambda_i(f(x_i)-f(\overline{x}))|\\
&\leq& |\sum_{i=1}^n \lambda_i(f(x_i)-f(\overline{x}))|,
\end{eqnarray*}
and therefore
$$||T||\leq \hbox{sup}\{|\sum_{i=1}^n\lambda_i(a_i-a)|: a_i,a\in [-1,1]\}.$$
Let $\phi:[-1,1]^{n+1}\longrightarrow \mathbb{R}$ be defined by 
$\phi(a_1,\ldots,a_n,a)=|\sum_{i=1}^n\lambda_i(a_i-a)|$. It is easy to see that $\phi$  
is continuous.
Moreover, since it is defined on a compact set, $\phi$  
assumes its maximum value at  $(b_1,\ldots,b_n,b)\in [-1,1]^{n+1}$. Then
$$||T||\leq |\sum_{i=1}^n\lambda_i(b_i-b)|.$$
On the other hand, consider $f$ in the unit ball of  $C(K)$ such that $f(x_i)=b_i$, for $i\leq n$, and $f(\overline{x})=b$.  
Recalling that $f_0(x_0)=1$, we get
$$|T(f)(x_0)|=|\sum_{i=1}^n \lambda_i(f(x_i)-f(\overline{x}))|=
|\sum_{i=1}^n\lambda_i(b_i-b)|\geq ||T||.$$ 
It follows that
\begin{eqnarray*}
||T||=||T(f)||= \hbox{sup}\{|\sum_{i=1}^n\lambda_i(a_i-a)|: a_i,a\in [-1,1]\}. 
\end{eqnarray*}
Consequently,  if  $n_0<n_1$ and  $\lambda_1, \ldots \lambda_{n_1}$ is a sequence of real numbers,  
we have 
$$\|\sum_{n=1}^{n_0} \lambda_nT_n\|\leq \|\sum_{n=1}^{n_1} \lambda_nT_n\|,$$
and hence the operators $T_n$ define a Schauder basis for $\overline{[T_i: n \geq 1]}$.
 
It remains to prove that each non-null element of that space is not a weak multiplier. 
 Let $T=\sum_{n=1}^\infty \alpha_nT_n$ be a non-null operator on $C(K)$  
and  let $n_0$ be the lowest natural number such that $\alpha_{n_0}\neq 0$.
Take a pairwise disjoint sequence $(e_n)_{n\in \N}$ in $C(K)$ 
such that $e_n(\overline{x})=0$ and
$$e_n(x_i)=\left\{ \begin{array}{ll}
			\frac{1}{\alpha_{n_0}}, & \mbox{se $i=n+n_0$}\\
			0, & \mbox{otherwise}
			\end{array}
			\right. $$
It is clear that $e_n(x_n)=0$, since $n\neq n+n_0$.  Moreover,  $T(e_n)(x_n)=1$ for all $n\in \N$.
Indeed, $T_{n_0}(e_n)(x_n)=e_n(x_{n+n_0})=\frac{1}{\alpha{n_0}}$, and
if $m\neq n_0$ we have $T_m(e_n)(x_n)=e_n(x_{n+m})=0$. It follows that $T\in NM(C(K))$.
\hfill{$\blacksquare$}

\vspace{2mm}

We will deal  now the general case. We will omit the proof of the following lemma to be quite simple.

\begin{lema} \label{lema:vaziolineavel} Let $(x_n)$ be a sequence in a compact space $K$ that does not converge. Then there exists  $f\in C(K)$  such that the sets
$$
A=\{ n\in\N : f(x_n)=0\} \hbox{ and } B=\{ n\in\N : f(x_n)=1\} 
$$ are both infinite.
\end{lema}

\begin{teor}\label{teor:vaziolineavel} 
The set of non weak multipliers on $C(K)$ is either empty or spaceable.
\end{teor} 

\paragraph{Proof:} If $K$ has a non-trivial  convergent sequence, the result follow from Theorem~\ref{teor:convergente}. Hence we may assume that $K$ does not contain a convergent sequence.

Suppose that there exists $T:C(K)\longrightarrow C(K)$ such that $T$
 is not a weak multiplier. Let $(e_n)_{n\in\N}$ be a bounded pairwise disjoint sequence in $C(K)$  and let $(x_n)$ be a sequence in $K$ such that $e_n(x_n)=0$ and $|T(e_n)(x_n)|>\varepsilon$,  for every $n\in \N$. 

Since $K$ does not contain  convergent sequences, using Lemma~\ref{lema:vaziolineavel}, we take  $f_1\in C(K)$  such that the sets
$$
A_1=\{ n\in\N : f_1(x_n)=0\} \hbox{ and } B_1=\{ n\in\N : f_1(x_n)=1\} 
$$ are both infinite. Assuming  $A_j, B_j$ and $f_j$  defined, again by the lemma, we take  $A_{j+1}, B_{j+1}$ and $f_{j+1}$ such that the sets 
$$
A_{j+1}=\{ n\in A_j : f_{j+1} (x_n)=0\} \hbox{ and  } B_{j+1}=\{ n\in A_j: f_{j+1}(x_n)=1\}
$$
are both infinite. 

For each $ j \geq 1 $, let $ T_j: C(K) \longrightarrow C(K) $ the operator defined by $ T_j (g) = f_j \cdot T(g) $ and consider  $S= \hbox{span}\{ T_j : j\geq 1\}$. 

Take  $\sum_{j=1}^k\alpha_j T_j\in S$, with $\alpha_k\not=0$.  Then for each  $n\in B_k$ we have 
$$
|\sum_{j=1}^k\alpha_j T_j(e_n)(x_n)| = |\alpha_k | |T_k(e_n) (x_n)|  > |\alpha_k| \varepsilon,
$$
since $B_k\subset A_j$, for all $j<k.$ This shows that the set $\{T_j: j \geq 1 \}$ is linearly independent and that, since each $B_k $  is infinite, each nonzero element of $S$ is not a weak multipliers.  So the set of non weak multipliers on $C(K)$ is lineable. 
By theorems~\ref{teor:closed} and~\ref{teor:equivalent} it is spaceable.
\hfill{$\blacksquare$}

\

\begin{ex} \rm If $C(K)$ has a proper subspace  $X$ isomorphic to it, then  $NM(C(K))$ is spaceable. Indeed, If $C(K)$ has a convergent sequence,  $NM(C(K))$ is spaceable by Theorem~\ref{teor:convergente}. Otherwise, by Theorem~\ref{teor:isosobre}, an isomorphism from  $C(K)$ onto  $X$ is not a  weak multipliers. So, by \ref{teor:vaziolineavel},    $NM(C(K))$ is spaceable. In particular,  $NM(C(\beta\N))=\ell_\infty$ is  spaceable, since is isomorphic to its square. Note that $\beta\N$  has no convergent sequence.
\end{ex}

As mentioned previously,  in \cite{Fe}  Ferenczi provides an example of  a Banach space $X$ on which there is an operator $J$ such that every operator on $X$ has the
form $\mu J+W$, where $\mu\in\mathbb{R}$ and $W=\lambda I+S$  with   $S$ strictly singular and  $\lambda\in\mathbb{R}$.  Note that the previous results show that this does not occur in the context    of  $C(K)$ spaces.   However, replacing scalar numbers by continuous functions we obtain something similar:

\begin{teor}\label{teor:almostfew} 
Let $K$ be a compact topological space such that every operator on $C(K)$
is a weak multiplication. Let $L=K\times\{0,1\}$. There exists an operator $J$ on $C(L)$
which is not weak multiplier such that every operator on $C(L)$ has the form $gJ+W$, where 
$g\in C(L)$ and $W$ is a weak multiplication. 
\end{teor}

\paragraph{Proof:} Let $K$ and $L$ be as in the hypothesis. 
For  $f\in C(L)$ we define $f_1, f_2 \in C(K)$ as  $f_1(x)=f(x,0)$ and $f_2(x)=f(x,1)$. We identify  $f\in C(L)$    with $(f_1,f_2)\in C(K)^2$. 
Define $J:C(L)\longrightarrow C(L)$ as $$J(f,g)=(g,f).$$ 

Let  $T$ be an operator on $C(L)$.  For each $i\in\{1,2\}$ we define  $\pi_i:C(L)\longrightarrow C(K)$ by $\pi_i(f)=f_i$ 
 and
  $\sigma_i:C(K)\longrightarrow C(L)$ by  
$\sigma_1(g)=(g,0)$ and $\sigma_2(g)=(0,g)$. 
Then, for each pair $(i,j)\in\{1,2\}^2$  we define an operator on $C(K)$ as 
$$T_{ij}=\pi_j\circ T\circ\sigma_i$$ 
and an operator on $C(L)$ as 
$$\tilde{T}_{ij}=\sigma_j\circ T_{ij}\circ\pi_i.$$ 
 
\begin{afi} $T=\sum_{i,j\in\{1,2\}} \tilde{T}_{i,j}.$ 
\end{afi} 
 
Let $f$ be an element of $C(L)$ which corresponds to $(f_1,f_2)$ in $C(K)^2$ 
(i.e., $\pi_i(f)=f_i$, for $i\in\{1,2\})$. 
We have \begin{eqnarray*} T(f)&=&(\pi_1(T(f_1,f_2)),\pi_2(T(f_1,f_2)))\\
&=& 
(\pi_1(T(f_1,0)+T(0,f_2)),\pi_2(T(f_1,0)+T(0,f_2)))\\
&=&(\pi_1(T(f_1,0))+\pi_1(T(0,f_2)),\pi_2(T(f_1,0))+\pi_2(T(0,f_2)))\\
&=& (\pi_1\circ T\circ\sigma_1(f_1)+\pi_1\circ T\circ\sigma_2(f_2), 
   \pi_2\circ T\circ\sigma_1(f_1)+\pi_2\circ T\circ\sigma_2(f_2) \\
&=& (T_{11}(f_1)+T_{21}(f_2),T_{12}(f_1)+T_{22}(f_2))\\
&=& (T_{11}(f_1),0)+(T_{21}(f_2),0)+(0,T_{12}(f_1))+(0,T_{22}(f_2))\\
&=& \sigma_1\circ T_{11}\circ\pi_1(f)+\sigma_1\circ T_{21}\circ\pi_2(f)+ 
\sigma_2\circ T_{12}\circ\pi_1(f)+\sigma_2\circ T_{22}\circ\pi_2(f)\\
&=&  \tilde{T}_{11}(f)+\tilde{T}_{21}(f)+\tilde{T}_{12}(f)+\tilde{T}_{22}(f), 
\end{eqnarray*}
which proves the claim. 
 
\vspace{2mm} 
 
To conclude the proof it is enough to show that each operator $\tilde{T}_{ij}$ 
has the form described in the theorem.  
 
\begin{afi} For any $i\in\{1,2\}$ there exists $\tilde{g}\in C(L)$ and a weakly compact operator 
$\tilde{S}$ on $C(L)$ such that $\tilde{T}_{ii}=\tilde{g}I_{C(L)}+\tilde{S}$. 
\end{afi} 
 
We may assume that $i=1$, since the proof is analogous for $\tilde{T}_{22}$. 
 
Since $T_{11}$ is an operator on $C(K)$ and $C(K)$ has few operators, there exist $g\in C(K)$ 
and an weakly compact operator $S$ on $C(K)$ such that $T_{11}=gI+S$.  
Take $\tilde{S}=\sigma_1\circ S\circ \pi_1$, i.e., $\tilde{S}(f_1,f_2)=(S(f_1),0)$. 
The operator $\tilde{S}$ is weakly compact, since composition of weakly compact operators 
with any operator is also weakly compact. 
 
Take $\tilde{g}=(g,0)\in C(L)$. For $f=(f_1,f_2)$ we have 
$$\tilde{T}_{11}(f)=(T_{11}(f_1),0)=(gf_1+S(f_1),0)=(gf_1,0)+(S(f_1),0)=\tilde{g}f+\tilde{S}(f).$$

\begin{afi} For any $i,j\in\{1,2\}$ such that $i\neq j$, there exist $\tilde{g}\in C(L)$  
and a weakly compact operator $\tilde{S}$ on $C(L)$ such that  
$\tilde{T}_{ij}=\tilde{g}J+\tilde{S}$. 
\end{afi} 
 
We proceed analogously to the above claim to prove that $\tilde{T}_{12}=\tilde{g}J+\tilde{S}$, for 
some $\tilde{g}\in C(L)$ and $\tilde{S}$ weakly compact on $C(L)$. 
 
Write $T_{12}=gI+S$. Take $\tilde{S}=\sigma_2\circ S\circ\pi_1$ and $\tilde{g}=(0,g)$. 
We have, for $f=(f_1,f_2)$,  
$$\tilde{T}_{12}(f)=(0,T_{12}(f_1))=(0,gf_1+S(f_1))=\tilde{g}\cdot (0,f_1)+\tilde{S}(f)= 
\tilde{g}J(f)+\tilde{S}(f).$$ 
 
This completes the proof of the claim and the theorem.  
\hfill{$\blacksquare$} 
 
\vspace{2mm} 

We notice that the operator $J$ in Theorem~\ref{teor:almostfew} can be replaced by the operator 
$J(f,g)=(-g,f)$, which has the property $J^2=-I$.

\bibliographystyle{amsplain}

\end{document}